\documentclass[12pt]{article}

\usepackage[margin=1in]{geometry}
\usepackage[T1]{fontenc}
\usepackage[utf8]{inputenc}
\usepackage{amsmath,amssymb,amsthm}
\usepackage{booktabs}
\usepackage{csquotes}
\usepackage[hidelinks]{hyperref}
\usepackage[style=apa,backend=biber,sorting=nyt,sortcites=true]{biblatex}
\usepackage[most]{tcolorbox}

\DeclareLanguageMapping{american}{american-apa}
\usepackage[
    backend=biber,
    style=apa
]{biblatex}

\newtcolorbox{activitybox}[1]{
    enhanced,
    breakable,
    colback=black!5,
    colframe=black!50,
    colbacktitle=black!12,
    coltitle=black,
    fonttitle=\bfseries,
    title={#1},
    boxrule=0.6pt,
    arc=1mm,
    left=8pt,
    right=8pt,
    top=8pt,
    bottom=8pt,
    before skip=12pt,
    after skip=12pt
}
\title{Scaffold--Write--Debug: A Three-Stage Approach to Proof Writing in Introductory Real Analysis}
\author{Chamila Malagoda Gamage\\
Department of Mathematics, University of Florida\\
\texttt{cgamage@ufl.edu}}
\date{}

\begin{document}

\maketitle

\begin{abstract}
This article describes Scaffold--Write--Debug, a three-stage instructional approach used in an introductory real analysis course to support students learning mathematical proof. The approach combines guided proof completion, independent proof writing, and the analysis and correction of faulty proofs. The goal is to help students understand proof structure, develop confidence in writing proofs, and become more careful readers of mathematical arguments. The paper presents classroom examples, instructor observations, and student responses, and discusses how the approach can be incorporated into a proof-based course.
\end{abstract}


\section{Introduction}

Introductory Real Analysis is often one of the first university mathematics courses in which students are expected to move from calculation to formal proof. In calculus, students may be asked to compute a limit, derivative, or integral and then check whether the final answer is reasonable. In real analysis, they are asked to work from definitions, choose a proof strategy, justify each step, and communicate the argument in a form that another reader can follow. This change can be difficult even for students who were successful in computational mathematics courses.

The challenge is not only understanding the mathematics. A student may understand a theorem but still not know how to start the proof. Another may have the right idea but leave gaps, misuse a definition, or introduce notation without explanation. Students may also follow a completed proof but struggle to write a similar one or decide whether an argument is correct. Proof writing therefore requires students to understand definitions, recognize the structure of a statement, choose useful results, build a clear argument, and check that each step follows logically.

A common approach is for the instructor to show complete proofs and then ask students to write new proofs on their own. Both are useful, but the jump between them can be too large. A finished proof shows what a correct argument looks like, but it does not always show how the proof was developed or why certain choices were made. On the other hand, an independent proof problem asks students to make many decisions at once. They must understand the statement, choose a method, recall definitions, organize the argument, and write it clearly. When they get stuck, they may not know exactly what is causing the difficulty.

To address this gap, the course used a three-stage approach called Scaffold–Write–Debug. Students first worked with partially completed proofs that made the main structure visible. They then used the same ideas to write related proofs independently. Afterward, they examined faulty proofs, identified where the reasoning failed, and rewrote the argument correctly. These stages were repeated across different topics, giving students regular practice in following, constructing, checking, and revising proofs.

This paper describes how the sequence was used in an introductory real analysis course. It includes examples from course notes and activities, together with instructor observations and student responses.

\section{Background and Related Literature}

\subsection{The Transition to Proof in Undergraduate Mathematics}

Learning to write proofs is different from solving familiar computational problems. In a proof problem, students must first understand the statement and decide which definitions, theorems, or methods may be useful. They must also recognize that examples, patterns, or diagrams may suggest an idea but do not prove it in general.

Research on the transition to proof has documented several sources of difficulty. \textcite{moore1994} identified difficulties related to conceptual understanding, mathematical language and notation, and getting started on a proof. These difficulties are closely connected. A student who does not understand a definition well may not recognize how it can be used in an argument. A student who understands the concept informally may still have difficulty expressing the idea with appropriate quantifiers or notation. A student who can follow an instructor's proof may not know which first step to choose when the structure is no longer supplied.

\textcite{weber2001} argued that proof construction requires strategic knowledge in addition to knowledge of definitions and theorems. Students need to know more than what mathematical facts are available. They also need ways to search for a proof, recognize a useful representation, connect the hypotheses to the desired conclusion, and decide when an approach should be changed. This distinction is important for introductory real analysis because students may memorize definitions of convergence, continuity, boundedness, or compactness while still having difficulty organizing those definitions into an argument.

The logical form of a statement also matters. \textcite{selden1995} showed the importance of unpacking mathematical statements so that their underlying logical structure becomes visible. In real analysis, a theorem may contain nested quantifiers, implications, or several conditions that must be coordinated. Students can lose track of which quantities are fixed, which may be chosen, and which depend on earlier choices. A proof can fail even when the student remembers the correct definition if the logical order of those choices is not preserved.

These studies show that students need both guidance and independent practice. They need help understanding proof structure, opportunities to make their own choices, and practice writing, checking, and revising arguments. The Scaffold--Write--Debug sequence addresses this need by moving students through different forms of participation in proof rather than expecting one activity to develop every part of proof competence.

\subsection{Scaffolding, Worked Examples, and Fading Support}

Scaffolding gives students temporary support for a task they are not yet ready to complete alone. As \textcite{wood1976} explain, this support helps learners manage difficult parts of a problem. In proof writing, it may include prompts, partially completed arguments, questions about definitions, or hints about what must be shown.

The support should gradually be reduced as students become more confident. Otherwise, students may learn only to complete a given structure rather than create one themselves. \textcite{atkinson2003} found that gradually removing worked-out steps, together with self-explanation prompts, can help students move toward independent problem solving. In proof writing, this means asking students not only to fill in missing steps but also to explain why those steps are valid.

Scaffolding can also take place through classroom discussion. \textcite{blanton2003} showed that instructor questions can help students compare ideas, explain their reasoning, and monitor their understanding. Similarly, \textcite{selden2018} used proof frameworks to make the logical structure of a theorem and its proof more visible. This can help students understand how definitions and assumptions guide the argument.

The skeletal proofs used in this course followed these ideas. Some provided the beginning and ending of a proof, while others included definitions or short prompts. As students gained experience, less support was provided. Students were also asked questions such as \enquote{Why does this follow?}, \enquote{Which definition is being used?}, and \enquote{Where is the hypothesis needed?} The goal was not only to complete the proof, but to understand how each step contributes to the argument.

\subsection{Independent Proof Writing and Practice}

Scaffolding is meant to lead students toward independent proof writing. In the write stage, students interpret a new statement, choose an approach, organize the argument, and explain their reasoning clearly.

Proof writing is both mathematical and communicative. A proof must be correct, but it must also be understandable to a reader. \textcite{lew2019} showed that students may not automatically understand common writing conventions, such as introducing variables, stating assumptions, and showing logical connections. For this reason, students should be given clear expectations for what a well-written proof should include.

Reading and explaining proofs can also support writing. \textcite{hodds2014} found that self-explanation improved proof comprehension, while \textcite{yee2018} showed the value of constructing, critiquing, and revising arguments. In the Scaffold--Write--Debug approach, students first work with support, then write independently, and later examine and correct faulty proofs. This allows them to practice both creating proofs and understanding what makes an argument clear and valid.

\subsection{Debugging, Proof Validation, and Learning from Errors}
The term \emph{debugging} is used here for checking, explaining, and correcting faulty proofs. Students are asked to locate the error, explain why it is wrong, and rewrite the proof correctly. Simply stating that a proof is incorrect is not enough.

Proof validation can be difficult for students. \textcite{selden2003} found that students may focus on the appearance of a proof rather than whether it actually proves the statement. Similarly, \textcite{alcock2005} showed that checking whether each line is true is not enough. Students must also check that each step follows logically from the previous work.

Debugging can also support proof writing. \textcite{powers2010} found that regular proof-validation activities can help students improve their own proofs. Work on incorrect examples also suggests that faulty solutions can support learning when students have enough background knowledge and are asked to explain the error carefully \parencite{grosse2007}. \textcite{pi2024} further showed that discussing, evaluating, and rewriting proofs can help students understand both correctness and clear mathematical writing.
\\

These components have been studied in different settings and have each been shown to support proof learning. Scaffolding helps students understand proof structure, independent writing gives them practice making their own choices, and validation and revision help them identify and correct weaknesses in an argument. \textcite{kirsten2023} also examined validation activities that occurred during undergraduate proof construction, including reviewing, correcting errors, expressing
doubts, and improving an argument. Their findings support the idea that validation is not only an activity performed after a proof is finished. It can be part of the construction process
itself. Building on this work, the present paper combines all three components into one repeated process intended to provide a practical and sustainable approach to learning proof writing in introductory real analysis.

\section{Course Context}

The course was an introductory real analysis course taught at a large public university in the United States. It had 25 students, most of whom were juniors or seniors majoring in mathematics. A few students had taken logic or other proof-based courses and had some experience with proof writing, although many were not confident in their abilities. For most students, this was their first course with a strong and regular emphasis on writing formal proofs.

The course covered familiar calculus topics, including limits, continuity, and derivatives, but approached them through abstract definitions, theorems, and proofs. It was taught during a regular 15-week semester and met three times per week for 50 minutes. Students also completed weekly homework, in-class exams, and class presentations.

\section{The Scaffold--Write--Debug Instructional Design}

The Scaffold--Write--Debug approach was used throughout the course to give students structured and repeated practice with proof writing. The three stages were connected to the same topic whenever possible. Students first worked with a partially completed proof, then wrote a related proof independently, and finally examined and corrected a faulty proof. The amount of support was reduced as students became more familiar with the definitions and proof techniques.

\subsection{The Scaffold Stage}

In the scaffold stage, students were given a skeletal proof with selected steps removed. The main structure of the argument was visible, but students had to complete important mathematical or logical steps. Some activities provided the beginning and end of the proof, while others included definitions, short hints, or questions about what should be shown next.

Students were also asked to explain their choices. For example, they might identify the definition being used, explain where a hypothesis was needed, or justify why one statement followed from another. This helped prevent the activity from becoming a simple fill-in-the-blank exercise. As the course continued, fewer steps and hints were provided.\\

\begin{activitybox}{Scaffold Activity: Proof by Induction}

\textbf{Exercise.}
Prove that for each natural number \(n\geq 1\),
\[
    \sum_{i=1}^{n} i
    =
    1+2+3+\cdots+n
    =
    \frac{n(n+1)}{2}.
\]

\bigskip

\textbf{Proof by induction}

\medskip

Let \(p(n)\) be the statement
\[
    \sum_{i=1}^{n} i=\frac{n(n+1)}{2}.
\]

\bigskip

\textbf{Step 1: Base case}

Verify the formula for \(n=1\).

\medskip

Left-hand side:
\[
    \sum_{i=1}^{1}i=\underline{\hspace{2cm}}
\]

Right-hand side:
\[
    \frac{1(1+1)}{2}=\underline{\hspace{2cm}}
\]

\medskip

\textbf{Is the left-hand side equal to the right-hand side?}
\[
    \underline{\hspace{3cm}}
\]

\bigskip

\textbf{Step 2: Induction hypothesis}

Assume that \(p(k)\) is true for some arbitrary integer \(k\geq 1\). That is,

\vspace{2cm}

\bigskip

\textbf{Step 3: Induction step}

We now prove that \(p(k+1)\) is true. That is,

\vspace{2cm}

Begin with
\[
    \sum_{i=1}^{k+1}i
    =
    \underline{\hspace{4cm}}.
\]

Using the induction hypothesis,
\[
    \sum_{i=1}^{k+1}i
    =
    \frac{k(k+1)}{2}+(k+1).
\]

Simplify the right-hand side.

\vspace{4cm}

\textbf{Conclusion.}
If \(p(k)\) is true, then \(p(k+1)\) is true.

\medskip

Therefore, by the principle of mathematical induction, \(p(n)\) is true for every
\(n\in\mathbb{N}\).

\end{activitybox}
\newpage
\subsection{The Write Stage}

In the write stage, students were asked to prove a related statement without a provided framework. They had to decide how to begin, which definitions or results to use, and how to organize the argument. These problems were usually connected to ideas students had already encountered in the scaffold stage, but they required students to make the main proof-writing decisions independently.

\begin{activitybox}{Write Activity: Proof by Induction}

\textbf{Exercise.}
Prove that for each natural number \(n\geq 1\),
\[
    \sum_{i=1}^{n} i^2
    =
    1^2+2^2+3^2+\cdots+n^2
    =
    \frac{n(n+1)(2n+1)}{6}.
\]

Write a complete proof by mathematical induction. Clearly include the base case, induction hypothesis, induction step, and final conclusion.

\end{activitybox}

\subsection{The Debug Stage}

In the debug stage, students were given a proof containing an error, gap, or unclear step. They were asked to locate the problem, explain why the reasoning was not valid, and rewrite the proof correctly. The errors were based on common difficulties and mistakes students make.

Early debugging activities usually contained one main error. Later activities sometimes included several related issues, such as an incorrect use of a definition, an unjustified inequality, or a conclusion that did not follow from the preceding argument. The correction step was important because students were expected not only to recognize a problem but also to show how the proof could be improved.\\

\begin{activitybox}{Debug Activity: An \(\varepsilon\)--\(\delta\) Proof}

\textbf{Claim.}
Prove that
\[
    \lim_{x\to 2}x^2=4.
\]

A student attempts to prove the claim using the \(\varepsilon\)--\(\delta\) definition.

\medskip

\textbf{Student's proof attempt.}

Let \(\varepsilon>0\). Choose
\[
    \delta=\frac{\varepsilon}{|x+2|}.
\]
Suppose that
\[
    0<|x-2|<\delta.
\]
Then
\[
    |x^2-4|
    =
    |x-2||x+2|
    <
    \delta|x+2|.
\]
Since
\[
    \delta=\frac{\varepsilon}{|x+2|},
\]
we obtain
\[
    |x^2-4|<\varepsilon.
\]
Therefore,
\[
    \lim_{x\to 2}x^2=4.
\]

\medskip

\textbf{Your tasks.}

\begin{enumerate}
    \item Identify the mistake in the student's proof.
    \item Explain why the choice of \(\delta\) is not valid.
    \item Rewrite the proof correctly using the
    \(\varepsilon\)--\(\delta\) definition.
\end{enumerate}

\end{activitybox}

\subsection{Using the Three Stages Together}

The three stages were used as a repeated cycle rather than as isolated activities. A scaffolded proof introduced the structure of an argument, an independent problem required students to apply the same ideas on their own, and a debugging activity helped them examine common errors more carefully. Not every topic required all three stages, and the activities were adjusted according to the difficulty of the material and students' previous experience. Overall, the sequence provided a gradual movement from guided work to independent proof writing and critical revision.

\section{Student Responses and Instructor Observations}
\label{sec:student_responses}

During the first week of the course, students completed a short survey about their previous experience with proofs, their confidence in reading and writing proofs, and the aspects of proof writing they found most difficult. The responses showed considerable variation in prior experience. Some students had taken a proof-based course before, while others had little formal experience with proof writing. Many students felt more comfortable reading and following an existing proof than constructing one from scratch. Commonly identified difficulties included ``knowing how to start or outline the argument,'' ``choosing a strategy,'' using definitions precisely, handling quantifiers and logical structure, and managing anxiety or time pressure. One student summarized a common goal simply: ``I want to get better at writing them!'' These responses confirmed my expectation that students would benefit from an explicit and sustained intervention in proof writing rather than being expected to develop the skill on their own.

As the semester progressed, students appeared less hesitant about beginning a proof and more willing to attempt an argument before knowing every step. They also became more familiar with common proof methods and more likely to begin by identifying the assumptions, the desired conclusion, and the relevant definitions. In their written work, I observed clearer organization, more careful use of definitions, and fewer unexplained logical steps. Students also became more willing to revise an argument after finding a gap or receiving feedback. These changes were gradual, but they suggested that proof writing was becoming a process that students could approach systematically rather than an unfamiliar task requiring an immediate clever idea.

A separate survey was given at the end of the semester. Most students who responded reported feeling more or much more confident in writing rigorous proofs than they had at the beginning of the course. They generally described their ability to structure a complete proof as moderate or strong. Their comments especially emphasized the value of having a clear starting point. One student wrote, ``The scaffolding gives you somewhere meaningful to start with any proof even if you don't really know where you're going, which is always helpful.'' Another identified the greatest improvement as ``getting started on the proofs and understanding the structure with which to write them.'' Students also reported that the scaffolded stage reduced the anxiety associated with facing a blank page and helped them recognize the general form of different proof techniques.

The individual writing and debugging stages supported different parts of the learning process. Writing proofs independently gave students an opportunity to determine whether they could apply the structure without guidance. The faulty-proof activities encouraged them to examine each step more critically and to consider whether a conclusion was fully justified. As one student explained, ``It helps me scrutinize my own proof better so that I know I'm not making logical jumps.'' Another reported being able to ``spot mistakes more easily and understand what logical flaws there were.'' Several students also indicated that the overall structure made proof writing less intimidating and helped them feel better prepared to write proofs independently.

The responses were not uniformly positive about every stage. Some students continued to struggle with choosing a method, identifying a useful intermediate step, or remembering a particular argument or ``trick.'' The debugging stage was also challenging for some students. One student wrote, ``It was a good exercise to think about possible mistakes, but I think it might have been good to have more practice with them.'' This suggests that error analysis may require repeated exposure before students become comfortable identifying subtle logical problems. Nevertheless, the overall feedback supports the value of combining scaffolding, independent practice, and debugging rather than relying on only one form of proof-writing instruction.

These survey responses are used descriptively and should not be interpreted as showing that the three-stage method alone caused the observed improvement. Students' development was also supported by regular instruction, homework, discussion, feedback, and continued practice throughout the semester. However, the surveys and classroom observations together suggest that the method provided students with a useful structure for approaching, writing, checking, and revising proofs.

\section{Practical Considerations and Conclusion}

The Scaffold--Write--Debug approach requires careful preparation. Instructors must decide how much of a proof to provide during the scaffold stage and when to remove that support. If too much is included, students may only complete the instructor's argument. If too little is included, the activity may not give enough guidance. The amount of scaffolding should therefore depend on the difficulty of the proof, students' prior experience, and the point in the semester.

Preparing useful faulty proofs also takes time. The errors should reflect natural mistakes that students may reasonably make, such as using a definition incorrectly, reversing an implication, choosing a quantity in the wrong order, leaving an important step unjustified, or failing to complete the conclusion. Artificial or obvious errors may turn the activity into a guessing exercise rather than careful mathematical reading. Faulty proofs should also contain a manageable number of problems so that students can locate, explain, and correct the main issue.

The approach does not need to be used in exactly the same way for every topic. Some proofs may require detailed scaffolding, while others may need only a short outline or reminder. Similarly, debugging activities are most useful after students have enough knowledge to recognize and repair the error. Instructors may also need to provide additional discussion when students can identify that a proof is wrong but cannot explain how to correct it.

Overall, the implementation of Scaffold--Write--Debug was successful in this introductory real analysis course. The sequence gave students a gradual way to move from following a proof to writing one independently and checking mathematical reasoning more carefully. Student responses suggested that the structure made proof writing less intimidating, helped them understand how to begin, and encouraged them to examine their own work more critically. Instructor observations also showed improvement in the organization, clarity, and completeness of students' proofs over the semester.

The approach does not remove all difficulties associated with proof writing, and students still need regular practice, feedback, and time to develop mathematical judgment. However, combining scaffolding, independent writing, and debugging provides a practical and sustainable way to support that development. More importantly, it presents proof writing not as a single skill that students either have or do not have, but as a process of understanding, constructing, checking, and revising mathematical arguments.

\printbibliography

\end{document}